\documentclass{amsproc}
\usepackage{graphicx}
\usepackage{epstopdf}
\usepackage{amsmath,amsthm,amscd,amssymb}
\usepackage{latexsym}
\usepackage[colorlinks,citecolor=red,pagebackref,hypertexnames=false]{hyperref}
\usepackage{geometry}                
\numberwithin{equation}{section}

\theoremstyle{plain}
\newtheorem{theorem}{Theorem}[section]
\newtheorem{lemma}[theorem]{Lemma}
\newtheorem{corollary}[theorem]{Corollary}
\newtheorem{proposition}[theorem]{Proposition}

\theoremstyle{definition}
\newtheorem{definition}[theorem]{Definition}

\newtheorem{problem}[theorem]{Problem}

\newtheorem{case[theorem]}{Case}

\theoremstyle{remark}
\newtheorem{remark}[theorem]{Remark}

\numberwithin{equation}{section}

\newcommand{\abs}[1]{\lvert#1\rvert}

\begin{document}

\title{\parbox{14cm}{\centering{Intersections of sets and Fourier analysis}}}


\author{Suresh Eswarathasan, Alex Iosevich, and Krystal Taylor}

\date{today}

\address{Institute des Hautes \'Etudes Scientifiques, 35 route des Chartres, Bures-sur-Yvette, France F-91440} 
\email{suresh@ihes.fr}
\address{Department of Mathematics and Statistics, McGill University, 805 Rue Sherbrooke Ouest, Montr\'eal, Canada H3A 0B9}
\email{suresh@math.mcgill.ca}
\address{Institute for Mathematics and its Applications, College of Science and Engineering, 207 Church Street SE, 306 Lind Hall, Minneapolis, MN U.S.A. 55455}
\email{krystaltaylormath@gmail.com}
\address{Department of Mathematics, University of Rochester, Rochester, NY}
\email{iosevich@math.rochester.edu}

\thanks{The first listed author was supported by a CRM-ISM Postdoctoral Fellowship and McGill University during the first part of the writing of this article, with the second part written while as a resident at the Institute des Hautes \'Etudes Scientifiques.  The work of the second listed author was partially supported by the NSF Grant DMS10-45404. The third listed author was supported by the Technion Technical Institute during the first part of the writing of this article, with second part written while at the Institute for Mathematics and its Applications.}

\begin{abstract} A classical theorem due to Mattila (see \cite{Mat84}; see also \cite{M95}, Chapter 13) says that if $A,B \subset {\Bbb R}^d$ of Hausdorff dimension $s_A, s_B$ respectively with $s_A+s_B \ge d$, $s_B>\frac{d+1}{2}$, and $dim_{{\mathcal H}}(A \times B)=s_A+s_B\ge d$, then 
$$ dim_{{\mathcal H}}(A \cap (z+B)) \leq s_A+s_B-d$$ for almost every $z \in {\Bbb R}^d$, in the sense of Lebesgue measure. 

In this paper, we replace the Hausdorff dimension on the left hand side of the first inequality above by the lower Minkowski dimension and replace the Lebesgue measure of the set of translates by a Hausdorff measure on a set of sufficiently large dimension.  Interesting arithmetic issues arise in the consideration of sharpness examples. These results are partly motivated by those in \cite{EIT11} and \cite{IJL10} where in the former the classical regular value theorem from differential geometry was investigated in a fractal setting, and in the latter discrete incidence theory is explored from an analytic standpoint.  

\end{abstract} 
\maketitle

\tableofcontents

\section{Introduction}

\vskip.125in 

A series of results due to Mattila (see \cite{Mat84}, \cite{Mat85}, \cite{Mat87}; see also \cite{M95}, Chapter 13) give lower and upper bounds on the Hausdorff dimension of the intersection of subsets of the Euclidean space of a given Hausdorff dimension.

\begin{theorem} \label{mattilalower} Let $A$ and $B$ be Borel subsets of $ {\Bbb R}^d$, $d\geq 2$, of Hausdorff dimension $s_A$ and $s_B$ respectively such that $s_A+s_B > d$ and $s_B>\frac{d+1}{2}$. Furthermore, assume that the respective Hausdorff measures of these sets are positive.  Then for almost every $g \in O(d)$, the group of orthogonal $d$ by $d$ matrices, 
\begin{equation}\label{mattilaconclusion} {\mathcal L}^d\left(\left\{z \in {\Bbb R}^d: dim_{{\mathcal H}}(A \cap (z-gB)) \ge s_A+s_B-d \right\}\right)>0.\end{equation}
\end{theorem} 

This means that for almost every rotation $g$ there is a set of $z$s of positive Lebesgue measure for which the Hausdorff dimension of $A \cap (z-gB)$ is at least $s_A+s_B-d$.  It is known that in a more restrictive setting,  if either $A$ or $B$ is Salem, then the assumption that one of the dimensions is at least $(d+1)/2$ is not necessary for $d\geq 2$ \cite{M95}.  In general, the necessity of this condition is not known.  In the case that $d=1$, setting both $A$ and $B$ equal to the middle-$\frac{1}{3}$ Cantor set demonstrates that \eqref{mattilaconclusion} fails in 1-dimension  \cite{M95}.  The converse does not in general hold, but the following result gives a partial description. 

\begin{theorem} \label{mattilaupper} Let $A$ and $B$ be Borel subsets of $ {\Bbb R}^d$, $d\geq 2$, of Hausdorff dimension $s_A$ and $s_B$ respectively
which satisfy
\begin{equation} \label{fair} dim_{{\mathcal H}}(A \times B)=s_A+s_B\ge d. \end{equation} 

Then \cite{M95}
\begin{equation} \label{upper} dim_{{\mathcal H}}(A \cap (z-B)) \leq s_A+s_B-d \end{equation} for almost every $z \in {\Bbb R}^d$ in the sense of Lebesgue measure. 
\end{theorem} 

This tells us that if (\ref{fair}) holds, then the Hausdorff dimension of $A \cap (z-gB)$ is at most $s_A+s_B-d$ for $g \in O(d)$ and almost every $z \in {\Bbb R}^d$. 
\vskip.125in 

We produce an example in remark \eqref{Federeradaptionkt} which illustrates that the assumption that \eqref{fair} holds in the previous theorem is a necessary condition.   
\vskip.125in 

A more general question, described in \cite{M95} and the references contained therein is the following. 

\begin{problem} \label{theproblem} To understand the Hausdorff dimension of $A \cap T(B)$, where $A,B$ are subsets of ${\Bbb R}^d$ of suitable Hausdorff dimension and $T$ ranges contained a suitable set of transformations of ${\Bbb R}^d$. \end{problem}

Before we give a detailed description of the goals of this paper, we wish to illustrate a simple motivating point by considering $A \cap (x-B)$, where $A,B \subset {\Bbb R}^d$. In order for the intersection to be non-empty, $x$ must be an element of the sum set $A+B$. If $A$ and $B$ are both sets of a given Hausdorff dimension $<d$, the Hausdorff dimension of $A+B$ is also quite often $<d$ and this naturally leads us to consider translates $x$ belonging to a set of a given Hausdorff dimension and exploring the thresholds for which the natural inequalities involving the dimension of $A \cap (x-B)$ hold. This simple point of view also indicates that the arithmetic properties of $A$ and $B$ play an important role. 
\vskip.125in 

An example of two sets $A$ and $B$, of Hausdorff dimension $s_A$ and $s_B$ respectively, such that the Hausdorff dimension of $A \cap (x-B)$ is ''generically" $s_A+s_B-d$ is easily constructed by taking $A$ and $B$ to be smooth surfaces in ${\Bbb R}^d$. A simple example in the non-integer case is obtained by considering 
$$ A=\{r \omega: \omega \in S^{d-1}; \ r \in U \},$$ where $U$ is an Ahlfors-David regular set of Hausdorff dimension $s_U$. It is not difficult to check that the Hausdorff dimension of $A$ is $d-1+s_U$. It is also straightforward to verify that the Hausdorff dimension of $A$ and every line that intersects $A$ is at most $s_U=d-1+s_U+1-d$. Modifying this construction yields examples of this type for arbitrary $s_A, s_B>\frac{d+1}{2}$, $s_A+s_B\geq d$. 
\vskip.125in

We now describe in some detail the goals of this paper: 

\vskip.125in 

\begin{itemize} 

\item Under structural assumptions on the sets $A$ and $B$, with $T_zB=B+z$, to prove that the set of translates $z$ for which the lower Minkowski dimension of $A \cap T_zB$ is larger than $dim_{{\mathcal H}}(A)+\dim_{{\mathcal H}}(B)-d$ does not only have Lebesgue measure $0$ but also a small Hausdorff dimension. This would be an analog of Theorem \ref{mattilaupper} above where finer information on the exceptional set of translates and replacing the Hausdorff dimension on the left hand side with lower Minkowski dimension is obtained at the expense of additional assumptions on the set 
$B$. 

\vskip.125in 

\item Without any additional assumptions on $A$ and $B$, beyond Ahlfors-David regularity, to replace the Hausdorff dimension by the lower Minkowski dimension in Theorem \ref{mattilaupper}. 

\vskip.125in 

\item To obtain the same type of results for $TB=gB+z$, where $g \in O(d)$ and $z \in {\Bbb R}^d$. We shall see that for almost every $g \in O(d)$, the set of translates $z$ for which lower Minkowski dimension of the set $A \cap TB$ is greater than $dim_{{\mathcal H}}(A)+dim_{{\mathcal H}}(B)-d$ has small Hausdorff dimension. 

\vskip.125in 

\end{itemize} 

\vskip.125in 



The main results of this paper are described in Section \ref{mainresults} below. The remainder of the paper is dedicated to proofs and remarks. 

\subsection{Notation}\label{notation} The following notation shall be used throughout: 

\vskip.125in 

\begin{itemize} 
\item $X \lesssim Y$ means that there exists $C>0$ which is independent of $X$ and $Y$ such that $X \leq CY$. 

\vskip.125in 

\item $X \lessapprox Y$ with the controlling parameter $R$ means that given $\epsilon>0$ there exists $C_{\epsilon}>0$ such that $X \leq C_{\epsilon}R^{\epsilon}Y$. 

\vskip.125in 

\item Let $B(x, \delta)$ be equal to the ball centered at a vector $x\in \mathbb{R}^d$ of radius $\delta$.  

\vskip.125in 

\item Given $A \subset {\Bbb R}^d$, let $A^{\epsilon}$ be equal to the open $\epsilon$-neighborhood of $A$. 

\vskip.125in 

\item Given $A \subset {\Bbb R}^d$  non-empty, let $s_A$ be equal to the Hausdorff dimension of $A$ and $\mu_A$ shall denote a probability measure on $A$. When $A$ is assumed to be Ahlfors-David regular, $\mu_A$ shall denote the restriction of the $s_A$-dimensional Hausdorff measure to $A$. 

\vskip.125in 

\item Given a compactly supported measure $\mu$ on ${\Bbb R}^d$, let $I^s(\mu)$ be equal to the $s$-energy integral given by $\int \int {|x-y|}^{-s} d\mu(x) d\mu(y)$.  By elementary properties of the Fourier transform, this expression is equivalent to $\int |\widehat{\mu}(\xi)|^2 |\xi|^{s-d}d\xi$.  Observe that if $0<\alpha<s$ and 
$\mu(B(x,\delta)) \leq C \delta^s$ for every $\delta>0$,
then $I^{\alpha}(\mu)\lesssim 1$ (see \cite{Fal86}, pg. 208 and \cite{Falc86}, section 6.2).   

\vskip.125in 

\item Given $A\subset \mathbb{R}^d$, let $M(A)$ be equal to the set of Radon measures $\mu$ with compact support such that the support of $\mu$ is contained in $A$ and $0< \mu(A) < \infty$.  
\end{itemize}

\vskip.125in

\section{Main results of this paper} 
\label{mainresults} 

\vskip.125in 

We shall primarily work with Ahlfors-David regular sets defined as follows. 

\begin{definition}\label{ADreg} We say that a Borel set $E \subset {\Bbb R}^d$ is Ahlfors-David regular if there exists $C>0$ such that for all $x \in E$
\begin{equation}\label{adreg} C^{-1} \delta^{s_E} \leq \mu(B(x, \delta)) \leq C \delta^{s_E},\end{equation} where $s_E$ is the Hausdorff dimension of $E$, $\mu$ is the Hausdorff measure restricted to $E$ and $B(x,\delta)$ is the ball of radius $0<\delta<diam(E)$ centered at $x$. 
\end{definition} 

\subsection{Intersections of translated, rotated and dilated sets} 
\label{generalstuff} 

\vskip.125in 

We begin with the following variant of Theorem \ref{mattilaupper}, where translation by $x \in {\Bbb R}^d$ is replaced by translation by $s(x)$, a local diffeomorphism and the Hausdorff dimension on the left hand side is replaced by the lower Minkowski dimension at the expense of assuming that the sets being intersected are Ahlfors-David regular. 

\begin{theorem} \label{eitmattilaupper}  Suppose that $A$ and $B$ are compact, Ahlfors-David regular, and Borel subsets of $ {\Bbb R}^d$, $d \ge 2$, of Hausdorff dimension $s_A$, $s_B$ respectively which satisfy $s_A + s_B \geq d$.

Now, let $\lambda(x)$ denote the lower Minkowski dimension of $A \cap (s(x)-B)$ 
where $s$ is a local $C^{\infty}$ diffeomorphism. Let $N(x,\epsilon)$ be equal to the minimum number of open $\epsilon$-balls needed to cover $A \cap (s(x)-B)$. 

Then for any smooth compactly supported function $\psi$
and every $\epsilon>0$, there exists a constant $C$, which is independent of $\epsilon$ and depends only on $\psi$ and $s$, such that
\begin{equation} \label{eitmattilaupperestquantitative} \int N(x,\epsilon) \psi(x)dx \leq C {(\epsilon^{-1})}^{s_A+s_B-d}, \end{equation} 
 from which it follows that 
\begin{equation}\label{eitmattilaupperest}\lambda(x) \leq s_A+s_B-d\end{equation} for almost every $x \in {\Bbb R}^d$. \end{theorem} 

\vskip.125in 

\begin{remark} In the case when $s(x) \equiv x$, Theorem \ref{eitmattilaupper} can be deduced from results in \cite{Falc94}. We include the result because our proof sets up the arguments in the remainder of the paper and due to its somewhat greater applicability. \end{remark} 

\begin{remark} The assumption that $B$ is Ahlfors-David regular in the statement of Theorem \ref{eitmattilaupper} can be eliminated at the cost of using the Minkowski dimension of $B$ on the right hand side of the inequality instead of the Hausdorff dimension of $B$. The authors are grateful to Brendan Murphy for making this observation.   \end{remark}


\begin{remark}\label{Federeradaptionkt} We wish to address the extent to which Theorem \eqref{eitmattilaupper} is sharp.  The following example illustrates that the dimensional inequality \eqref{eitmattilaupperest} fails when we remove the assumption that the sets $A$ and $B$ are Ahlfors-David regular. In fact, we show something even stronger;  That is, if 
$$dim(A \times B) > dim(A) + dim(B),$$ which implies that neither $A$ nor $B$ is Ahlfors-David regular, then \eqref{eitmattilaupperest} fails even if the lower Minkowski dimension of the intersection set is replaced with the Hausdorff dimension of the intersection set.  \\

Let $A,B \subset \mathbb{R}^2$ be defined as follows:
$$A=[0,2] \times Y,$$ and $$B=X \times [0,2]$$ where $X,Y \subset [0,1]$ are Borel sets which 
satisfy 
 \begin{equation}\label{keytodemo}dim(X \times Y) > dim(X) + dim(Y).\end{equation}
The existence of such sets is discussed in  \cite{M95} and \cite{Fe69}. Federer constructs examples of such sets in  (\cite{Fe69}, (2.10.29)) which are also Borel subsets of $[0,1]$.  \\
Observe that $$A\cap B = ([0,2] \cap X) \times (Y\cap [0,2]) =X \times Y,$$
and $$A\cap (B+ (u,v)) = ([0,2] \cap (X+u)) \times (Y\cap [v,2+v]) =(X+u)\times Y,$$ whenever $(u,v) \in [0,1]\times [-1,0].$
Since $X$ and $Y$ are Borel, it follows by a result in  \cite{M95} that 
$$dim(A) = 1+ dim(Y)$$ and $$dim(B) = 1+ dim(X).$$
Combining these observations, we conclude that
\begin{align}
& \label{demokt}dim(A\cap (B+ (u,v)) ) = dim((X+u) \times Y)=dim(X\times Y) \\
&  > dim(X) + dim(Y) = dim(A) + dim(B) -2,
\end{align} whenever $(u,v) \in [0,1]\times [-1,0].$  In other words, \eqref{eitmattilaupperest} fails on a set of positive measure.  
\end{remark}

Combining Theorem \ref{eitmattilaupper} with Theorem \ref{mattilalower}, we deduce that the lower Minkowski dimension and the Hausdorff dimension of the intersection of an Ahlfors-David regular set with a rotated copy of a Borel set quite frequently coincide. 

\begin{corollary} \label{adregular} Suppose that $A$ and $B$ are compact, Ahlfors-David regular, and Borel subsets of $ {\Bbb R}^d$, $d \ge 2$, of Hausdorff dimension $s_A$, $s_B$ respectively which satisfy $s_A + s_B \geq d$.  Also assume that 
 $s_B>\frac{d+1}{2}$. Then for almost every $g \in O(d)$, 
\begin{equation} \label{adregularconclusion} {\mathcal L}^d \left\{z \in {\Bbb R}^d: dim_{{\mathcal H}}(A \cap (z-gB))=\underline{dim}_{{\mathcal M}}(A \cap (z-gB)) \right\}>0. \end{equation} 
\end{corollary} 

\vskip.125in 

\begin{remark} It is reasonable to conjecture that under the assumptions of Corollary \ref{adregular}, $A \cap (z-gB)$ is Ahlfors-David regular, but this does not follow from the equality of the Hausdorff and lower Minkowski dimensions. This can be seen by taking a Cantor construction and changing the dissection ratio at each stage. The second listed author is grateful to Pertti Mattila for pointing this construction in the context of Ahlfors-David regularity. \end{remark} 

\vskip.125in 

If we are willing to rotate $B$ before translating it, we discover that the exceptional set, which was found to have Lebesgue measure zero in Theorem \ref{eitmattilaupper} above, has a small Hausdorff dimension. 

\begin{theorem} \label{rotated}  Suppose that $A$ and $B$ are compact, Ahlfors-David regular, and Borel subsets of $ {\Bbb R}^d$, $d \ge 2$, of Hausdorff dimension $s_A$, $s_B$ respectively which satisfy $s_A + s_B \geq d$. Let $\mu$ be a compactly supported probability measure such that $I^{\alpha}(\mu)<\infty$, for some $0<\alpha\le d$ satisfying  
$$ \alpha+s_A>2(d-\kappa)$$ where $\kappa=\min\{ \frac{d-1}{2}, s_B\}$.   

Let 
$\lambda_{g}(x)$ denote the lower Minkowski dimension of $A \cap (x-gB)$, where $g \in O(d)$, the orthogonal group of transformations on $\mathbb{R}^d$. Let $N(x,g,\epsilon)$ be equal to the minimum number of open $\epsilon$-balls needed to cover $A \cap (x-gB)$. Then there exists $C>0$ such that for every $\epsilon>0$
\begin{equation} \label{rotatedconclusionquantitative} \int \int N(x,g, \epsilon) d\theta(g) d\mu(x) \leq C \sqrt{I^{\alpha}(\mu)I^{2(d-\kappa)-\alpha}(\mu_A)} \cdot {(\epsilon^{-1})}^{s_A+s_B-d}, \end{equation} where $d\theta(g)$ denotes normalized Haar measure on $O(d)$.

It follows that 
\begin{equation} \label{rotatedconclusion}\lambda_g (x)  \leq s_A+s_B-d \end{equation} 
almost everywhere with respect to the probability measure $d\theta(g)d\mu(x)$.\\

Furthermore, 
\begin{equation} \label{rotatedsizzle} dim_{{\mathcal H}}\left(\left\{x:  \int \lambda_g(x) d\theta(g)>s_A+s_B-d \right\}\right)\le d+1-s_A. \end{equation} 

\end{theorem} 

\vskip.125in 

We are also able to obtain a good upper bound on the Hausdorff dimension of the exceptional set if we put additional structural assumptions on one of the sets being intersected.

\vskip.125in 

\begin{definition} We say that a Borel and compact set $B \subset {\Bbb R}^d$ of Hausdorff dimension $s_B$ satisfies the hyperplane size condition of order 
$h$, for some $s_B>h>0$, if there exists a Borel measure $\mu_B$ supported in $B$ such that 
$$ \mu_B\left(H^{\delta}_{\omega}\right) \leq C{\delta}^{s_B-h},$$ where $H_{\omega}=\{x \in {\Bbb R}^d: x \cdot \omega=0\}$. 
\end{definition} 

\begin{remark} Note that if $\mu_B$ is a Frostman measure, then the hyperplane size condition with $h=d-1$ always holds. This is because the intersection of $B$ with a hyperplane can be decomposed into $\approx \delta^{-(d-1)}$ $\delta$-cells, and the measure of  each cell is $\leq C\delta^{s_B}$ by the Frostman property. One should think of $h$ as an upper bound on the dimension of the intersection of $B$ with a $(d-1)$-dimensional hyperplane. \end{remark} 

\begin{theorem} \label{dilation}   Suppose that $A$ and $B$ are compact, Ahlfors-David regular, and Borel subsets of $ {\Bbb R}^d$, $d \ge 2$, of Hausdorff dimension $s_A$, $s_B$ respectively which satisfy $s_A + s_B \geq d$.

 Let $\mu$ be a compactly supported probability measure with $I^{\alpha}(\mu)<\infty$. Furthermore, suppose that $B$ satisfies the hyperplane size condition of order $h<s_B$ and 
$$ \frac{\alpha+s_A}{2}>d-(s_B-h).$$ 

Let $\lambda_t(x)$ denote the lower Minkowski dimension of $A \cap (x-tB)$. Let $N(x,t,\epsilon)$ be equal to the minimum number of open $\epsilon$-balls needed to cover $A \cap (x-tB)$. Then 
\begin{equation} \label{dilationconclusionquantitative} \int \int_1^2 N(x,t, \epsilon) dt d\mu(x) \leq C {(\epsilon^{-1})}^{s_A+s_B-d}. \end{equation} 

It follows that 
\begin{equation} \label{dilationconclusion} \lambda_t(x) \leq s_A+s_B-d \end{equation} almost everywhere with respect to the probability measure $dt \, d\mu(x)$ and $t\in[1,2]$.

Finally, for almost every $t \in [1,2]$,
\begin{equation} \label{dilationsizzle} dim_{{\mathcal H}}\left(\left\{x:\int_1^2 \lambda_t(x)dt >s_A+s_B-d\right \}\right)\le 2(d-(s_B-h))-s_A. \end{equation}

\end{theorem} 

\vskip.25in 

\vskip.25in 

\section{Proof of Theorem \ref{eitmattilaupper}}
Let $\mu_A$ and $\mu_B$ denote the restrictions of the $s_A$ and $s_B$-dimensional Hausdorff measures to $A$ and $B$ respectively which are normalized so that $\int d\mu_A=\int d\mu_B=1$, and take $\psi \in C^{\infty}_0(\mathbb{R}^d)$.  
In order to prove \eqref{eitmattilaupperestquantitative}, we obtain upper and lower bounds on the following quantity 
\begin{equation}\label{special}
\int \mu_A({\{A\cap (s(x)-B)\}^{\epsilon}})\psi(x) dx.
\end{equation} 
We first obtain a lower bound, and we utilize the following lemma:

\begin{lemma} \label{lemmauppersemi}
Let $A$, $B$, and $s$ be as in Theorem \eqref{eitmattilaupper}.  
For a fixed $\epsilon>0$, set $g(x)=N(x,\epsilon)$, where $N(x,\epsilon)$ is equal to the minimum number of open $\epsilon$-balls needed to cover $A \cap (s(x)+B)$.
Then $g(x)$ is an upper semi-continuous function on $\mathbb{R}^d$.  As a result, $N(x,\epsilon)$ is measurable on $\mathbb{R}^d$.
\end{lemma}

We delay the proof of this lemma until the appendix for purposes of fluidity.  Multiplying each side of the following equation 
$$N(x,\epsilon)\epsilon^{s_A} \lesssim \mu_A (\{A\cap (s(x)-B)\}^{\epsilon} )$$ by $\psi(x)$ and integrating in $x$, which is allowed by Lemma \ref{lemmauppersemi}, we use the Ahlfors-David regularity of $A$ to conclude  that  \eqref{special} is bounded below by a dimensional constant times $$ \epsilon^{s_A}\int  N(x,\epsilon) \psi(x) dx.$$  

We next establish an upper bound on the expression \eqref{special}.  Observe that ${(A \cap (s(x)-B))}^{\epsilon} \subset \{y \in A^{\epsilon}: s(x)-y \in B^{\epsilon} \}$, and therefore \eqref{special} is bounded above by  
\begin{equation}\label{main2}
\int \mu_A(\{y \in A^{\epsilon}: s(x)-y \in B^{\epsilon} \})\psi(x)dx.
\end{equation}
Let $J_s$ denote the Jacobian of the change of variables $x \to s(x)$. Notice that

$$ \int \mu_A \{y \in A^{\epsilon}: s(x)-y \in B^{\epsilon} \} \psi(x) dx $$

$$ \approx \int \int \chi_{B^{\epsilon}}(s(x)-y) d\mu_A(y) \psi(x) dx .$$ 
Here, $\chi_{B^{\epsilon}}$ is some cutoff supported in a small neighborhood of $B^{\epsilon}$.  
We use properties of the Fourier transform to obtain the following bound:
$$ \left|\int \int \chi_{B^{\epsilon}}(s(x)-y) d\mu_A(y) \psi(x) dx\right| \le \int |\widehat{\chi_{B^{\epsilon}}}(\xi)|  \cdot |\widehat{\mu}_A(\xi)| \cdot 
| \left(J_{s^{-1}} \cdot (\psi \circ  s^{-1})\right)^{\widehat{}} (\xi)| d\xi.  $$ 

Let $\mathcal{L}^d(B^{\epsilon})$ denote the $d$-dimensional Lebesgue measure of the set $B^{\epsilon}$.  Since $\|\widehat{\mu}_A\|_{\infty}\le \int d \mu_A(x)=1$, $\psi\circ s^{-1}$ is a smooth and compactly supported function, and $\|\widehat{\chi_{B^{\epsilon}}}\|_{\infty}\le \|\chi_{B^{\epsilon}}\|_{1}$, we further bound the expression above by
\begin{equation}\label{expression} 
C  \cdot \mathcal{L}^d(B^{\epsilon})
\end{equation} 
where $C$ is independent of $\epsilon$. \\
Let $N(B,\epsilon)$ be equal to the minimum number of open  $\epsilon$-balls needed to cover $B$. 
Recall that the Hausdorff dimension of $B$ is equal to the Minkowski dimension of $B$. 
This follows, for instance, by the Ahlfors-David regularity of $B$.
Using the definition of the Minkowski dimension, we see that $\mathcal{L}^d(B^{\epsilon})\sim \epsilon^d N(B,\epsilon) \sim \epsilon^{d-s_B} $. 
Now the expression in \eqref{expression} is $\sim$ $c\epsilon^{d-s_B}.$ 

Comparing our upper and lower bounds, it follows that for any smooth compactly supported function $\psi$, any smooth diffeomorphism $s$ of $\mathbb{R}^d$,
and for every $\epsilon>0$, we have that
$$\int N(x,\epsilon) \psi(x)dx \lesssim {\epsilon}^{-(s_A+s_B-d)}.$$ This is precisely \eqref{eitmattilaupperestquantitative}.

The proof of \eqref{eitmattilaupperest} is a technical exercise in limit inferiors of sequences, and is delayed until the appendix for the purpose of continuity.

\vskip.125in 

\section{Proof of Theorem \ref{rotated}} 

\vskip.125in 
Let $\mu_A$ and $\mu_B$ denote the restrictions of the $s_A$ and $s_B$-dimensional Hausdorff measures to $A$ and $B$ respectively which are normalized so that $\int d\mu_A=\int d\mu_B=1$. In order to prove \eqref{rotatedconclusionquantitative}, we obtain upper and lower bounds on the following quantity
\begin{equation}\label{boundthis}  \int \int \mu_A (\{A\cap (x-gB)\}^{\epsilon} ) d\mu(x) d\theta(g) .\end{equation}
Let $N(x, g,\epsilon)$ be equal to the minimum number of open $\epsilon$-balls needed to cover $A\cap(x-gB)$.  
By the properties of $\mu_A$, which are a consequence of the Ahlfors-David regularity of $A$, it follows that
$$N(x,g,\epsilon)\epsilon^{s_A} \lesssim \mu_A (\{A\cap (x-gB)\}^{\epsilon} ).$$ 
Integrating each side of this equation in $x$ and $g$ with respect to the measures $d\mu$ and $d\theta$ respectively, we conclude that \eqref{boundthis} is bounded below by
\begin{equation}\label{lowerboundrotated}\int \int N(x,g,\epsilon)\epsilon^{s_A} d\mu(x) d\theta(g),\end{equation}
where the measurability of $N(x,g,\epsilon)$ follows from arguments similar to those used in the proof of Lemma \ref{lemmauppersemi}, which appears in the appendix, after working in local coordinates on $O(d)$.

Next, observe that $\{A\cap( x-gB)\}^{\epsilon} \subset \{y \in A^{\epsilon}: x-y \in gB^{\epsilon} \}$, and so \eqref{boundthis} is bounded above by 
$$ \int \int \mu_A \{y \in A^{\epsilon}: x-y \in gB^{\epsilon} \} d\mu(x) d\theta(g).$$  
This quantity is comparable to 

\begin{equation} \label{grh+} \int \int \int \chi_{B^{\epsilon}}(g(x-y)) d\mu(x) d\mu_A(y) d\theta(g)
\end{equation} 
$$=\int \left( \int \widehat{\chi_{B^{\epsilon}}}(g\xi) d\theta(g) \right) \overline{\widehat{\mu}(\xi)} \widehat{\mu}_A(\xi) d\xi.$$ 
We use the following lemma to bound the previous expression.
\begin{lemma} \label{propeller} 
With the notation above, 
$$  \left| \int \widehat{\chi_{B^{\epsilon}}}(g\xi) d\theta(g) \right| \lesssim \epsilon^{d-s_B}{(1+|\xi|)}^{-\kappa},$$
where $\kappa=\min\{\frac{d-1}{2}, s_B\}$.   
\end{lemma} 

Postponing the proof of Lemma  \ref{propeller} until the end of this section,  we can now bound the expression in (\ref{grh+}) above by
$$   \epsilon^{d-s_B} \int {|\xi|}^{-\kappa} |\widehat{\mu}(\xi)| \cdot |\widehat{\mu}_A(\xi)| d\xi.$$  Write $-\kappa =\frac{(\alpha-d)}{2}  + \left( -\kappa  -\frac{(\alpha-d)}{2} \right)$.  
By the Cauchy-Schwartz inequality, the previous expression is bounded by 
\begin{equation}\label{upperboundrotated}  \epsilon^{d-s_B} \sqrt{I^{\alpha}(\mu)I^{2(d-\kappa) -\alpha}(\mu_A)}\lesssim  \epsilon^{d-s_B} .\end{equation}
Indeed, $I^{\alpha}(\mu)$ is finite by assumption, and $I^{2(d-\kappa)-\alpha}(\mu_A)\lesssim 1$ follows by the hypothesis that 
$s_A + \alpha>2(d-\kappa)$.  \\

Combining the upper and lower bounds for \eqref{boundthis}, found in \eqref{upperboundrotated} and \eqref{lowerboundrotated} respectively, it follows that 
$$\epsilon^{s_A}\cdot \int \int N(x,\epsilon) d\mu(x) d\theta(g) \lesssim \sqrt{I^{\alpha}(\mu)I^{2(d-\kappa) -\alpha}(\mu_A)} \epsilon^{d-s_B}, $$ which is precisely \eqref{rotatedconclusionquantitative}.
\\

Next, observe that  \eqref{rotatedconclusion}  follows from \eqref{rotatedconclusionquantitative}, and the proof is almost identical to that of (\ref{eitmattilaupperest}) in Theorem \ref{eitmattilaupper}.  The only changes that need to be made are that $dx$ is replaced by $d\theta(g) d\mu(x)$, and the words ``positive Lebesgue measure" are replaced by ``positive measure with respect to $d\theta(g) d\mu(x)$".
To demonstrate \eqref{rotatedsizzle}, assume by way of contradiction that the exceptional set $\{ x :  \lambda_{g}(x)  > s_A +s_B -d \}$ has dimension larger than $d +1  -s_A$.  We are also assuming that $d+1-s_A< d$ as otherwise the claim holds trivially.  Choose $\alpha$ such that $$d+1-s_A < \alpha < \dim_{\mathcal{H}}\{ x : \int \lambda_{g}(x) d\theta(g) > s_A +s_B -d \}.$$  By Frostman's Lemma \cite{M95}, there exists a compactly supported probability measure $\mu$ with support contained in this exceptional set so that $I^{\alpha}(\mu)<\infty$.  We use the observation that \eqref{rotatedconclusion} holds for the measure $\mu$ to arrive at a contradiction. This concludes the proof of our theorem, but it remains to prove Lemma \ref{propeller}. 

\subsection{Proof of Lemma \ref{propeller}}
We have 
\begin{align*} 
&\int \widehat{\chi_{B^{\epsilon}}}(g\xi) d\theta(g)\\
&=\int \int e^{-2 \pi i g^{-1}x \cdot \xi} d\theta(g)\chi_{B^{\epsilon}}(x) dx\\
&=c\int \widehat{\sigma}(|x| \xi) \chi_{B^{\epsilon}}(x) dx,
\end{align*}
where $\sigma$ is the Lebesgue measure on $S^{d-1}$ and $c$ is a constant depending only on $B$. 
Using the well-known estimate of the Fourier Transform of the sphere, $\abs{  \widehat{\sigma}(\xi) }\lesssim (1+ |\xi|)^{-\frac{(d-1)}{2}}$, we bound the modulus of the expresion above by a constant times
\begin{equation}\label{propellermain} \int (1+(\abs{x}\abs{\xi}))^{-\frac{(d-1)}{2}} \chi_{B^{\epsilon}}(x) dx.\end{equation}
We will consider two cases:  First, we consider the case when $|\xi|\le 1$, and then we consider the case when $|\xi|>1$.  
When $|\xi|\le1$,  we bound the expression in \eqref{propellermain} by
$$\int \chi_{B^{\epsilon}}(x) dx.$$
Now, we approximate the $d$-dimensional Lebesgue measure of the set $B^{\epsilon}$ by the minimal number of $\epsilon-$balls needed to cover $B$ times the size of such a ball, $\epsilon^d$.  
Since $B$ is Ahlfors-David regular, it follows that the Hausdorff dimension of $B$ is equal to the lower Minkowski dimension of $B$, and so the $d$-dimensional Lebesgue measure of $B^{\epsilon}$ is approximately $\epsilon^{d-s_B}$.  This shows that the lemma holds when $\abs{\xi}\le 1$.\\

Next, we consider the case when $\abs{\xi}>1$, and we break the integral in \eqref{propellermain} over three regions: 
$\left\{x: \abs{x}<\frac{1}{\abs{\xi}}\right\}$, $\left\{x: \frac{1}{\abs{\xi}}<\abs{x}<1\right\}$, and $\left\{x: 1<\abs{x}\right\}$.
That is, for $\abs{\xi}>1$, we estimate \eqref{propellermain} from above by 
\begin{align*}
& \int_{\left\{x: \abs{x}<\frac{1}{\abs{\xi}}\right\}} \chi_{B^{\epsilon}}(x) dx 
\, +\int_{\left\{x: \frac{1}{\abs{\xi}}<\abs{x}<1\right\}} (\abs{x}\abs{\xi})^{-\frac{(d-1)}{2}} \chi_{B^{\epsilon}}(x) dx  \\
& + \int_{\left\{x: 1<\abs{x}\right\}}(\abs{x}\abs{\xi})^{-\frac{(d-1)}{2}} \chi_{B^{\epsilon}}(x) dx \\
& =: I + II + III.
\end{align*}
In each region, we will utilize the following elementary estimates:
\begin{proposition}\label{propellerprop}
Let $\rho: \mathbb{R}^d \rightarrow \mathbb{R}$ be a non-negative, smooth, and compactly supported function which is greater or equal to one on the ball of radius two centered at the origin.  Take $B$ to be an Alhfors-David regular set.  For $\epsilon>0$, let $\rho_{\epsilon}(x) = \frac{1}{\epsilon^d}\rho(\frac{x}{\epsilon})$. Then, for $\delta>0$, we have the following estimates:
\begin{align}
 \label{lemmabab2}\chi_{B^{\epsilon}}(x) & \lesssim \epsilon^{d-s_B}\mu_{B}*\rho_{\epsilon}(x), \\
\label{lemmababy2}   \int_{\{x: \abs{x}<\delta\}} \mu_{B}*\rho_{\epsilon}(x) dx & \lesssim \delta^{s_B},\\
\label{lemmababy3} \int \mu_{B}*\rho_{\epsilon}(x) dx & \lesssim 1.
\end{align}
\end{proposition}

These observations are standard and are left as an exercise for the reader.  Now, we use Proposition \ref{propellerprop} to estimate \eqref{propellermain} restricted to the first region:
\begin{align*}
I
&= \int_{\{x: \abs{x}<\frac{1}{\abs{\xi}}  \}} \chi_{B^{\epsilon}}(x) dx
\\
& \lesssim  \epsilon^{d-s_B}  \int_{\{x: \abs{x}<\frac{1}{\abs{\xi}} \}} \mu_{B}*\rho_{\epsilon}(x) dx
\\
& \lesssim \epsilon^{d-s_B} \abs{\xi}^{-s_B}.
\\
\end{align*}
Next, we estimate \eqref{propellermain} restricted to the second region:
\begin{align*}
II &=\int_{\{x: \frac{1}{\abs{\xi}}<\abs{x}<1\}} (\abs{x}\abs{\xi})^{-\frac{(d-1)}{2}} \chi_{B^{\epsilon}}(x) dx
\\
&\lesssim \epsilon^{d-s_B }\int_{\{x: \frac{1}{\abs{\xi}}<\abs{x}<1\}} (\abs{x}\abs{\xi})^{-\frac{(d-1)}{2}} \mu_{B}*\rho_{\epsilon}(x) dx
\\
&\sim   \epsilon^{d-s_B }\abs{\xi}^{-\frac{(d-1)}{2}} \sum_{j=0}^{\log_2(\frac{1}{\abs{\xi} }) -1} \int_{\{x: 2^{-(j+1)}<\abs{x}<2^{-j}\}}  \abs{x}^{-\frac{(d-1)}{2}} \mu_{B}*\rho_{\epsilon}(x) dx
\\
&\sim \lesssim \epsilon^{d-s_B }\abs{\xi}^{-\frac{(d-1)}{2}} \sum_{j=0}^{\log_2(\frac{1}{\abs{\xi} }) -1}  2^{j\left(\frac{(d-1)}{2}-s_B\right)} 
\\
&\lesssim \epsilon^{d-s_B}|\xi|^{-\kappa},\\
\end{align*}
where $\kappa=\min\{\frac{d-1}{2}, s_B\}$.
Finally, we use Propositon \ref{propellerprop} to estimate \eqref{propellermain} restricted to the third region:
\begin{align*}
III&=\int_{\{x: 1<\abs{x}\}} (\abs{x}\abs{\xi})^{-\frac{(d-1)}{2}} \chi_{B^{\epsilon}}(x) dx
\\
&\lesssim \epsilon^{d-s_B } \abs{\xi}^{-\frac{(d-1)}{2}}  \int_{\{x: \abs{x}>1\}} \abs{x}^{-\frac{(d-1)}{2}}   \mu_{B}*\rho_{\epsilon}(x)   dx
\\
&\lesssim \epsilon^{d-s_B } \abs{\xi}^{-\frac{(d-1)}{2}}  \int \mu_{B}*\rho_{\epsilon}(x)   dx
\\
&\lesssim \epsilon^{d-s_B } \abs{\xi}^{-\frac{(d-1)}{2}}.
\\
\end{align*}
This concludes the proof of Lemma \ref{propeller}. 
\vskip.125in

\vskip.125in 

\section{Proof of Theorem \ref{dilation}}

In order to prove \eqref{dilationconclusionquantitative}, we obtain upper and lower bounds on the following quantity:
 \begin{equation}\label{quantity4}\int \int_1^2 \mu_A({\{A\cap( x-tB)\}^{\epsilon}}) dt d\mu(x).\end{equation}
Let $N(x,t,\epsilon)$ be equal to the minimum number of open $\epsilon-$balls needed to cover $A\cap (x-tB)$.
Observe the following which holds by the Ahlfors-David regularity of the set $A$:
\begin{equation}\label{obs1kt} N(t,x,\epsilon) \epsilon^{s_A} \lesssim \mu_A({\{A\cap (x-tB)\}^{\epsilon}}).\end{equation}
Integrating each side of this equation with respect to the measure $dtd\mu(x)$, $t\in (1,2)$, we conclude that \eqref{quantity4} is bounded below by 
\begin{equation}\label{lower4} \epsilon^{s_A} \int \int_1^2 N(x,t,\epsilon) dt d\mu(x). \end{equation}

Next observe that  $\{A\cap (x-tB)\}^{\epsilon} \subset A^{\epsilon}\cap (x-tB^{\epsilon})$. 
 Hence, \eqref{quantity4} is bounded above by
\begin{equation}\label{main13}\int \int_1^2 \mu_A({A^{\epsilon}\cap(x-tB^{\epsilon})})dt d\mu(x) .\end{equation}

Let $\rho: \mathbb{R}^d \rightarrow \mathbb{R}$ be a non-negative smooth-bump function which is greater or equal to one on $B(\vec{0}, 1/2)$ and equal to zero outside of $B(\vec{0}, 1)$.  It follows that, 

\begin{equation}\label{obs4} \mu_A({A^{\epsilon}\cap(x-tB^{\epsilon}))} \lesssim \epsilon^{d-s_B}\int\mu_B * \rho_{4\epsilon}\left(\frac{x-y}{t}\right)d\mu_A(y).\end{equation}

The measurability of $N(x,t, \epsilon)$ follows, once again, from arguments similar to that of Lemma \ref{lemmauppersemi}. Integrating in $x$ with respect to the measure $\mu$ and in $t\in[1,2]$ we bound \eqref{quantity4} from above by the following expression:
\begin{equation}\label{compareEkt}\epsilon^{-(s_B-d)}\int \int \int \mu_B * \rho_{4\epsilon}\left(\frac{x-y}{t}\right) \psi(t) d\mu_A(y)dt d\mu(x),\end{equation}
where $\psi$ is a translated smooth bump function equal to one on $[1,2]$ and equal to zero outside of $[0.5, 2.5]$.

\begin{lemma}\label{dilationlemmakt} With the notation above, 
\begin{equation}\label{dilationkt}\iint \int \mu_B * \rho_{\epsilon}\left(\frac{x-y}{t}\right) \psi(t) d\mu_A(y)dt d\mu(x)\lesssim 1,\end{equation}
whenever $\frac{s_A+\alpha}{2}> d-(s_B-h)$. \end{lemma} 

Assuming Lemma \ref{dilationlemmakt} for the moment,  we have demonstrated that \eqref{quantity4} is bounded above by a constant times $\left(\frac{1}{\epsilon}\right)^{s_B-d}$ and is bounded below by a constant times the expression in \eqref{lower4}.  That is, 
$$ \epsilon^{s_A} \int \int_1^2 N(x,t,\epsilon) dt d\mu(x) \lesssim \epsilon^{-(s_B-d)},$$ and  \eqref{dilationconclusionquantitative} is proved. 

Once again, we observe that (\ref{dilationconclusion}) follows from  \eqref{dilationconclusionquantitative}, and the proof is almost identical to that of (\ref{eitmattilaupperest}) in Theorem \ref{eitmattilaupper}.  The only changes that need to be made are that $dx$ is replaced by $dt \, d\mu(x)$ for $t\in [1,2]$, the smooth cut off function $\psi$ is replaced with the constant function equal to one everywhere, and the words ``positive Lebesgue measure" are replaced by ``positive measure with respect to $dt \, d\mu(x)$".

To demonstrate \eqref{dilationsizzle}, we follow same method as in the previous proofs.  We assume once again, by way of contradiction, that the exceptional set $\{ x :  \lambda_t(x) > s_A +s_B -d \}$ has dimension larger than $2d-2s_B+2h-s_A$.  We are also assuming that $2d-2s_B+2h-s_A<d$ as otherwise the claim holds trivially.  Let $$2d-2s_B+2h-s_A< \alpha < \dim_{\mathcal{H}}(\{ x : \lambda_t(x)  > s_A +s_B -d \}).$$  By Frostman's Lemma \cite{M95}, there exists a compactly supported probability measure $\mu$ with support contained in this exceptional set so that $I^{\alpha}(\mu)<\infty$.  We use  \eqref{dilationconclusion} to arrive at a contradiction.

In order to finish the proof of the theorem, it remains to prove Lemma \ref{dilationlemmakt}.  Using elementary properties of the Fourier transform, we re-write the left-hand-side of \eqref{dilationkt} as 
\begin{equation}\label{ftkt}\int\int  \left (\mu_B*\rho_{\epsilon}\right)^{\widehat{}}(t\xi)\widehat{\mu_A}(\xi)\overline{\widehat{\mu}(\xi)} t^d\psi(t)dt d\xi.\end{equation}
The modulus of this expression is bounded above by a constant times
\begin{equation}\label{modkt}\int \abs{\widehat{\mu_A}(\xi)}\abs{\widehat{\mu}(\xi)} 
 \left| \int \left (\mu_B*\rho_{\epsilon}\right)^{\widehat{}}(t\xi)t^d\psi(t)dt \right| d\xi.\end{equation}
The following estimate is a key point towards establishing our lemma.

\begin{proposition}\label{hyperplaneprop} 
With the notation above, for all $s_B-h>\eta>0$ there exists $c_{\eta}$, which is independent of $\epsilon$, such that
\begin{equation}\label{hyperplaneequation}
 \left| \int \left (\mu_B*\rho_{\epsilon}\right)^{\widehat{}}(t\xi)t^d\psi(t)dt \right| \lesssim c_{\eta}\abs{\xi}^{-(s_B-h)+\eta}.
\end{equation}
\end{proposition}
Postponing the proof of Proposition \ref{hyperplaneprop} to the end of this section, we can now bound the expression in \eqref{modkt} above by a constant times
\begin{equation}\label{integral1kt}
 \int |\widehat{\mu_A}(\xi)|\cdot |\widehat{\mu}(\xi)| \cdot |\xi|^{-(s_B-h)+\eta} d\xi\end{equation}

Write $-(s_B-h) =\frac{\alpha-d}{2} - (s_B-h)- \frac{\alpha-d}{2}$.  It follows by the Cauchy-Schwartz inequality that \eqref{integral1kt} is bounded above by 
\begin{equation}\label{energykt}
\sqrt{ I^{\alpha}(\mu) I^{2(d-(s_B-h)+\eta)-\alpha}(\mu_A)}.\end{equation} 
By assumption, $I^{\alpha}(\mu)\lesssim 1$.  Also, observe that  $s_A + \alpha>2(d-(s_B-h))$  implies that $I^{2(d-(s_B-h)+\eta)-\alpha}(\mu_A)\lesssim 1$ for $\eta>0$ choosen sufficiently small.  We conclude that \eqref{energykt} is bounded by a positive constant which does not depend on $\epsilon$ whenever $\frac{s_A+\alpha}{2}> d-(s_B-h)$.  This completes the proof of Lemma \ref{dilationlemmakt} up to the proof of Proposition \ref{hyperplaneprop}.  We devote a subsection towards proving the estimate \eqref{hyperplaneequation}.

\subsection{Proof of Proposition \ref{hyperplaneprop}}

Consider 
\begin{equation}\label{tintkt}\int \left (\mu_B*\rho_{\epsilon}\right)^{\widehat{}}(t\xi) t^d\psi(t)dt.\end{equation}
Motivated by the presence of $\widehat{\rho_{\epsilon}}(t\xi)$, we consider the case when $\abs{\xi}<\frac{1}{\epsilon}$ and the case when $\abs{\xi}>\frac{1}{\epsilon}$ separately.  
\\

Consider the case when $\abs{\xi}<\frac{1}{\epsilon}$.
Set $\tilde{\psi} =t^d \psi$ and $\mu_B^{\epsilon}= \mu_B*\rho_{\epsilon}$. 
 Use the definition of the Fourier transform to re-write \eqref{tintkt} as
\begin{equation}\label{changekt}
\int\left (\mu_B^{\epsilon}\right)^{\widehat{}}(t\xi) \tilde{\psi}(t)dt =\int \widehat{\tilde{\psi}}(x\cdot \xi) \mu_B^{\epsilon}(x)dx. \end{equation}
Next, we use the rapid decay of $\tilde{\psi}$ to bound the modulus of this expression by a constant $C_N$ times
 \begin{equation}\label{moddecay}
\int  (1+ |x\cdot \xi|)^{-N} \mu_B^{\epsilon}(x)dx, \end{equation}
for any $N>1$.  We separate this integral over two regions, namely $\left\{x : \abs{x\cdot \xi} <1\right\}$ and $\left\{x : \abs{x\cdot \xi} >1\right\}$.
That is, the expression in \eqref{moddecay} can be estimated above by 
\begin{align*}
& \int_{\{x : \abs{x\cdot \xi} <1\}}   \mu_B^{\epsilon}(x)dx   + \int_{\{x : \abs{x\cdot \xi} >1\}}   |x\cdot \xi|^{-N} \mu_B^{\epsilon}(x)dx\\
&= I + II.
\end{align*}

We first consider \eqref{moddecay} restricted to the first region and write out the integrand explicitly:
\begin{align*}
I &= \int_{\{x : \abs{x\cdot \xi} <1\}}   \mu_B^{\epsilon}(x)dx \\
&= \frac{1}{\epsilon^d} \int \int_{\{x: |x\cdot \xi | \le 1\}} \rho\left( \frac{x-y}{\epsilon}\right)dx d\mu_B(y) \\
\end{align*}
Fix $y$, recall that $\rho$ is supported in $B(\vec{0},1)$, and break the integral in $x$ to consider $\{x: |x-y| < \epsilon\}$ and $\{ x: 1>|x-y| > \epsilon\}$.  
That is, for any $M>1$, we have the further decomposition
\begin{align*}
I &\lesssim \frac{1}{\epsilon^d} 
\int \int_{\{x: |x\cdot \xi | \le 1\}\cap \{x: |x-y| < \epsilon\}} dx d\mu_B(y)
+ \int \int_{\{x: |x\cdot \xi | \le 1\}\cap \{x: 1>|x-y| > \epsilon\}} \left| \frac{x-y}{\epsilon}\right|^{-M}dx d\mu_B(y)\\
&= I_a + I_b.\\
\end{align*}

Observe that
\begin{align*}
I_a
& =\frac{1}{\epsilon^d} \int \int_{\{x: |x\cdot \xi | \le 1\}\cap \{ x: |x-y| < \epsilon\}}dx d\mu_B(y) \\
&\lesssim \mu_B( \{y : |y\cdot \xi| \le 1 + \epsilon |\xi|\} ). \\
&\lesssim \mu_B\left( \left\{y : \left|y\cdot \frac{\xi}{|\xi|}\right| \le \frac{2}{|\xi|}\right \}\right),\\
\end{align*}
where the last line follows since we assumed that $|\xi| < \frac{1}{\epsilon}$.\\
Applying the hyperplane size condition of order $h$ on the set $B$, we conclude that
$$I_a \lesssim |\xi|^{-(s_B -h)}.$$

Similarly, we have that
\begin{align*}
I_b
&= \frac{1}{\epsilon^d} \int \int_{\{x: |x\cdot \xi | \le 1\}\cap \{ x: 1>|x-y| > \epsilon\}}\left|\frac{x-y}{\epsilon}  \right|^{-M} dx d\mu_B(y) \\
&\approx \epsilon^{M-d}\sum_{j=0}^{\log_2(\frac{1}{\epsilon})-1}\int \int_{\{x: |x\cdot \xi | \le 1\}\cap \{ x: 2^{-(j+1)}<|x-y|<2^{-j}\}} 2^{jM} dx d\mu_B(y) \\
&\lesssim \epsilon^{M-d} \sum_{j=0}^{\log_2(\frac{1}{\epsilon})-1} 2^{jM}  2^{-jd} \mu_B( \{y : \left|y\cdot \frac{\xi}{\abs{\xi}}\right| \le 2^{-j} + \epsilon \} ). \\
\end{align*}
Once again applying the hyperplane size condition of order $h$ on the set $B$, we may bound this quantity by
$$\epsilon^{M-d} \sum_{j=0}^{\log_2(\frac{1}{\epsilon})-1} 2^{j(M-d-(s_B-h))}\approx \epsilon^{s_B-h},$$
for $M$ sufficiently large. 
Reminding ourselves that $|\xi| < \frac{1}{\epsilon}$, we conclude 
$$I_b \lesssim |\xi|^{-(s_B-h)}.$$

Moving onto the second region in the case of $|\xi| < \frac{1}{\epsilon}$, we consider:
\begin{align*}
II&= \int_{\{x : \abs{x\cdot \xi} >1\}}    |x\cdot \xi|^{-N}  \mu_B^{\epsilon}(x)dx. \\
\end{align*} 
By breaking the integral further into the regions $\left\{x: 2^k <|x\cdot \xi|<2^{k+1}\right\}$, for $k\in \mathbb{N}$, we are able to show that 
$$II \lesssim |\xi|^{-(s_B-h)}$$
with nearly identical estimates to those used for bounding $I$.  We leave this as an exercise for the reader.\\

It remains to consider the case when $\abs{\xi}>\frac{1}{\epsilon}$.  
We will consider the following two subcases: $\abs{\xi}>\left(\frac{1}{\epsilon}\right)^{1+1/c}$ and  $\left(\frac{1}{\epsilon}\right)^{1+1/c} > \abs{\xi}>\frac{1}{\epsilon}$, where a positive lower bound on $c$ will be determined and shown to be independent of $\epsilon$.  
We obtain two estimates on the modulus of \eqref{tintkt} which will be used in either case. Re-visiting the estimates for $I_a$, $I_b$, and $I_c$ with the assumption that $\abs{\xi}>\frac{1}{\epsilon}$, we see that the modulus of the expression in \eqref{tintkt} is bounded by another dimensional-constant times $\epsilon^{s_B-h}$.  On the otherhand, using the rapid decay of $\widehat{\rho}$, we bound the modulus of the expression in \eqref{tintkt} by a dimensional-constant times $C_N\cdot(\epsilon|\xi|)^{-N}$ for any $N\geq 1$.  Indeed,
\begin{align*}
&\int |\widehat{\mu_B}(t\xi)| | \widehat{\rho}(\epsilon t\xi)| t^d\psi(t)dt\\
&\lesssim C_N \int (1+\epsilon t|\xi|)^{-N} t^d\psi(t)dt,\\
\end{align*}
for any $N\geq 1$.  
Because we are assuming that $\abs{\xi} >\frac{1}{\epsilon}$ and because $\psi(t)=0$ outside of $[0.5,2.5]$, we see that $\epsilon t|\xi| \geq \frac12\epsilon|\xi| $ on $[1,2]$, and so we have the upper bound on \eqref{tintkt} over the indicated region:
$$C_N (\epsilon|\xi|)^{-N}\int t^{d-N}\psi(t)dt\lesssim c_N (\epsilon|\xi|)^{-N}.$$

Set $N=(c+1)(s_B-h)\geq 1$, where an additional positive lower bound on $c$ will be choosen momentarily.  If $\abs{\xi}> \left(\frac{1}{\epsilon}\right)^{1+1/c}$, then one may verify that $(\epsilon \abs{\xi})^{-N} < \abs{\xi}^{-(s_B-h)},$ and so we bound \eqref{tintkt} by $c_N\cdot(\epsilon|\xi|)^{-N} \lesssim c_N \abs{\xi}^{-(s_B-h)}.$ Hence, Proposition \ref{hyperplaneprop} holds with constants independent of $\epsilon$ for $\abs{\xi}> \left(\frac{1}{\epsilon}\right)^{1+1/c}$.

If $\abs{\xi}>\left(\frac{1}{\epsilon}\right)^{1+1/c}$, then we bound \eqref{tintkt} by $\epsilon^{s_B-h}$, and we see that for this range of $\abs{\xi}$ it holds that
$$\epsilon^{s_B-h}\lesssim \abs{\xi}^{-(s_B-h) + \frac{s_B-h}{c+1}}.$$  Choosing $c$ sufficiently large so that $\frac{s_B-h}{c+1}<\eta$, we conclude our proof of Proposition \ref{hyperplaneprop}.

\section{Appendix}

\subsection{Proof of \eqref{eitmattilaupperest}}

We will use \eqref{eitmattilaupperestquantitative} to prove \eqref{eitmattilaupperest} using a proof by way of contradiction. Set $$C=\left\{x: \lambda(x) > s_A + s_B -d\right\},$$ and assume that the $d$-dimensional Lebesgue measure of $C$ is positive.  This implies that there exists a real number $N>0$ such that 
\begin{equation}\label{contradictionset}C_{N}=\left\{x: \lambda(x) > s_A + s_B -d+ \frac{2}{N}\right\}\end{equation} also has positive $d$-dimensional Lebesgue measure. 
We restrict our attention to $x\in C_{N}$.

Begin by re-writing the set $C_N$.  For $j \in \mathbb{N}$, define
\begin{equation}\label{firstdecomp}D_{N,J}=\left\{  x \in C_{N}: \lambda(x) -\frac{1}{N}<  \frac{\log(N(x,\epsilon))}{\log(\frac{1}{\epsilon})} \text{ for all }0<\epsilon\le 2^{-j}   \right\}.\end{equation} Observe that  $$C_N =\bigcup_{j=1}^{\infty}D_{N,j}.$$ 
To see this, fix $x\in C_{N}$ and  recall that, for $A\cap(s(x) -B)\neq \emptyset$,
\begin{align*}
\lambda(x)
& = \liminf_{\epsilon\downarrow0}\left(  \frac{\log(N(x,\epsilon))}{\log(\frac{1}{\epsilon})}\right) \\
&= \lim_{\delta\downarrow0}\left(\inf\left\{  \frac{\log(N(x,\epsilon))}{\log(\frac{1}{\epsilon})} : 0<\epsilon\le\delta \right \} \right).\\
\end{align*}

Certainly, $A\cap(s(x) -B)\neq \emptyset$ since  $x\in D_{N,j} \subset C_N$ implies that $\lambda(x)>s_A +s_B-d+\frac{2}{N}>0$.    \\
By the definition of the limit, there exists $j \in \mathbb{N}$ such that 
$$\lambda(x) -\frac{1}{N}<\inf\left\{  \frac{\log(N(x,\epsilon))}{\log(\frac{1}{\epsilon})} : 0<\epsilon \le 2^{-j} \right \}.$$ 
By the definition of the infimum, it follows that 
$$\lambda(x) -\frac{1}{N}<  \frac{\log(N(x,\epsilon))}{\log(\frac{1}{\epsilon})}$$ for all $0<\epsilon \le 2^{-j}$. 
This establishes that  $$C_N =\bigcup_{j=1}^{\infty}D_{N,j}.$$ 

Notice that $D_{N,j} \subset D_{N,j+1}$.  Recalling that $\mathcal{L}^d(C_N)>0$ and $C_N =\bigcup_{j=1}^{\infty}D_{N,j}$, it follows that 
$$\mathcal{L}^d\left(\bigcup_{j=1}^J D_{N,j}\right)=\mathcal{L}^d(D_{N,J})>0$$
for some $J$ sufficiently large.  

To summarize, we have found a set $ D_{N,J}\subset C_N$ of positive Lebesgue measure and a $J$ sufficiently large, such that $x\in D_{N,J}$  implies 
\begin{equation}\label{work} s_A+s_B-d +\frac{1}{N} < \lambda(x) -\frac{1}{N} < \frac{\log(N(x,\epsilon))}{\log(\frac{1}{\epsilon})}\end{equation} for all $\epsilon \in \left(0,2^{-J}\right]$.
It follows that $x\in D_{N,J}$ and $\epsilon\in (0,2^{-J}]$ implies that 
\begin{equation}\left(\frac{1}{\epsilon}\right)^{s_A+s_B-d +\frac{1}{N}} <N(x,\epsilon).\end{equation}

Define $\psi_J$ to be a smooth, compactly supported, and non-negative function such that $$\int_{D_{N,J}} \psi_J(x) dx =1.$$ This is possible as $D_{N,J}$ has positive Lebesgue measure. We then obtain that
\begin{equation}\label{knockoutma1}\int_{D_{N,J}} \left(\frac{1}{\epsilon}\right)^{(s_a+s_b-d+\frac{1}{N})} \psi_J(x)dx 
< \int_{D_{N,J}} N(x,\epsilon) \psi_J(x)dx.\end{equation} whenever $\epsilon \in \left(0,2^{-J}\right]$.  Since  $\int_{D_{N,J}} \psi_J(x) dx =1$, it follows that 
\begin{equation}\label{knockoutmama1}
\left(\frac{1}{\epsilon}\right)^{(s_A+s_B-d+\frac{1}{n})}
< \int_{D_{N,J}} N(x,\epsilon) \psi_J(x)dx.\end{equation}

Using \eqref{eitmattilaupperestquantitative} to bound the right-hand-side of this expression, we obtain
\begin{equation}\label{knockoutmamama1}
 \left(\frac{1}{\epsilon}\right)^{(s_A+s_B-d+\frac{1}{N})}
< \int_{D_{N,J}} N(x,\epsilon) \psi_J(x)dx < C'\left(\frac{1}{\epsilon}\right)^{(s_A+s_B-d)} \end{equation}  whenever $\epsilon \in \left(0,2^{-J}\right]$ and $C'>0$ is independent of $\epsilon$.  If we choose $\epsilon$  sufficiently small, then  \eqref{knockoutmamama1} cannot hold, and we arrive at a contradiction.  Therefore, it must hold that $$\lambda(x) \le s_A + s_B -d$$  for almost every $x \in \mathbb{R}^d$ with respect to Lebesgue measure.

\subsection{Proof of Lemma \ref{lemmauppersemi}}

In order to prove that $g(x)$ is upper semi-continous at a point $x_0$, we find $\delta>0$ so that $|x-x_0|<\delta$ implies that $g(x)\le g(x_0)$.   That is,  we fix $\epsilon>0$ and find a value of $\delta>0$ such that $|x-x_0|<\delta$ implies that
\begin{equation}\label{uppersemi}  N(x,\epsilon) \le N(x_0,\epsilon).\end{equation}
We consider the case when $N(x_0,\epsilon)=0$ and $N(x_0,\epsilon)\neq0$ separately.  In both cases, we will utilize the following proposition which will be proved at the end of this section:
\begin{proposition}\label{simple}
Let $X,U$ be non-empty sets in $\mathbb{R}^d$ such that $X\subset U$, $X$ is compact, and $U$ is open.  Then, there exists $\delta>0$ so that \begin{equation}\label{wiggle}X^{\delta}\subset U.\end{equation} 
\end{proposition}

We now turn to the proof of \eqref{uppersemi} in the case when $N(x_0,\epsilon)=0$.  Note that $N(x_0,\epsilon)=0$ occurs if and only if $A\cap(s(x_0)+B)= \emptyset$. In this case, $s(x_0)+B\subset \mathbb{R}^d\backslash A$.  Since $s(x_0)+B$ is closed and $ \mathbb{R}^d\backslash A$ is open, it follows by Proposition \ref{simple} that there exists $\lambda>0$ so that $\{s(x_0)+B\}^{\lambda}\subset  \mathbb{R}^d\backslash A$.
By the continuity of $s$, there exists $\delta>0$ so $|x-x_0|<\delta$ implies that $|s(x) -s(x_0)|< \lambda$, and so $s(x) + B \subset \{s(x_0)+B\}^{\lambda}\subset  \mathbb{R}^d\backslash A$.  
We have shown that there exists $\delta>0$ so that $|x-x_0|<\delta$ implies that $A\cap (s(x)+B)=\emptyset$, and \eqref{uppersemi} is established for this case.  \\

Now, we prove \eqref{uppersemi} in the case that $N(x_0,\epsilon)\neq 0$. For simplicity of notation, set $N=N(x_0,\epsilon)$.  Then there exist $N$ open $\epsilon-$balls, denote them as $B_{\epsilon}(c_1)\cdots B_{\epsilon}(c_N)$, so that 
\begin{equation}\label{cover} A\cap (s(x_0)+B) \subset \bigcup_{i=1}^N B_{\epsilon}(c_i).\end{equation}
We show that there exists $\delta>0$ so that $|x-x_0|<\delta$ implies that 
\begin{equation}\label{mainuppersemi}A\cap (s(x)+B) \subset \bigcup_{i=1}^N B_{\epsilon}(c_i).\end{equation}
This is enough to prove \eqref{uppersemi} because it would follow from \eqref{mainuppersemi} that $$N(x,\epsilon)\le N.$$ 

In order to establish \eqref{mainuppersemi}, we first apply Proposition \ref{simple} to \eqref{cover} to conclude that there exists $\lambda>0$ so that 
\begin{equation}\label{lemma2uppersemi} \{A\cap (s(x_0)+ B)\}^{\lambda}\subset \bigcup_{i=1}^{N}B_{\epsilon}(c_i).\end{equation}

We next show that there exists $\delta>0$ so that whenever $|x-x_0|<\delta$ then 
\begin{equation}\label{lemma1uppersemi}  A\cap(s(x)+B) \subset \{A\cap (s(x_0)+B)\}^{\lambda}.\end{equation}
  \\
Note that \eqref{mainuppersemi} would follow from  \eqref{lemma2uppersemi} and \eqref{lemma1uppersemi}.  In order to prove  \eqref{lemma1uppersemi}, we use the continuity of $s$ and the following proposition:  

\begin{proposition}\label{simpleuppersemi} Let $\lambda>0$ so that \eqref{lemma2uppersemi} holds.   For any non-empty and compact sets $A$ and $B$, there exists $\eta>0$ so that whenever $|x|<\eta$ then 
\begin{equation}  A\cap(x+B) \subset \{A\cap B\}^{\lambda}.\end{equation}
\end{proposition}  

Delaying the proof of Proposition \ref{simpleuppersemi} momentarily, we now deduce  \eqref{lemma1uppersemi} from  the proposition and the continuity of $s$. \\

Indeed, replacing $B$ in the statement of Proposition \ref{simpleuppersemi} with $s(x_0)+B$, it follows  that there exists $\eta>0$ such that $|y|<\eta$ implies that $A\cap (y+s(x_0)+B) \subset \{A\cap (s(x_0)+B)\}^{\lambda}.$ By the continuity of $s$, there exists $\delta>0$ so that $|x-x_0|<\delta$ implies that $|s(x) -s(x_0)|<\eta$.  That is, $|x-x_0|<\delta$ implies that $s(x)=y+s(x_0)$ for some $|y|<\eta$ and so $A\cap (s(x) + B) \subset \{A\cap (s(x_0) + B)\}^{\lambda}$.  This establishes \eqref{lemma1uppersemi} modulo the proofs of Proposition \ref{simpleuppersemi} and Proposition \ref{simple}.

\subsection{Proof of Proposition \ref{simple}}
We begin by covering each point in $X$ by a ball of radius $r_x$, $B(x,r_x)$, so that the ball of radius $2r_x$ is contained in $U$ (this is possible since $U$ is open).  
Next, we extract a finite subcovering from this open covering of $X$:  $$X\subset \bigcup_{i=1}^N B(x_i,r_i)\subset \bigcup_{i=1}^N B(x_i,2r_i)\subset U.$$ 
Let $\delta>0$ equal the minimum of the $r_i$ for $i=1,\cdots,N$. For this choice of $\delta$ it follows that 
$$X^{\delta}\subset U.$$  To see this, let $v\in X^{\delta}$.  Then, $|v-x|<\delta$ for some $x\in X$. It follows that $x\in B(x_i,r_i)\subset X$ for some $i=1,\cdots,n$, and so $v\in B(x_i,r_i+\delta)\subset B(x_i,2r_i) \subset U$.

\subsection{Proof of Proposition \ref{simpleuppersemi}}

We begin by creating a finite cover of the compact set $B$.  We will consider $b \in A\cap B$ and $b\in\left(\mathbb{R}^d\backslash A\right) \cap B$ separately. 
Fix $\lambda>0$, and for each $b\in A\cap B$, place a ball of radius $\frac{\lambda}{2}$ centered at b. 
For each  $b\in\left(\mathbb{R}^d\backslash A\right) \cap B$, cover $b$ by a ball centered at $b$ of radius $\gamma_b$ so that $B(b,2\gamma_b)\cap A =\emptyset$ (this is possible because $A$ is closed).  From this cover, extract a finite subcover of $B$: 
$$B\subset \left( \bigcup_{i=1}^N B(b_i,\lambda/2)\right) \bigcup \left(\bigcup_{i=1}^M B(\tilde{b_i},\gamma_i)  \right),$$ 
where $b_i \in A\cap B$ and each $\tilde{b_i}$ is at distance at least $2\gamma_i$ from $A$.  

Let $\delta=\min\{\frac{\lambda}{2},\gamma_1,\cdots,\gamma_M\}$.   Now $|x|<\delta$ guarantees that $$A\cap (x+B) \subset \{A\cap B\}^{\lambda}.$$  Indeed, let $v=x+b\in A$ where $b\in B$ and $|x|<\delta$.  Clearly, $b\notin B(\tilde{b_i},\gamma_i)$ for $i=1,\cdots,M$ (since otherwise, $|\tilde{b_i}-v|\le |\tilde{b_i}-b| +|b-v| < \gamma_i + \delta \le 2\gamma_i$ which contradicts $B(\tilde{b_i},2\gamma_i)\cap A =\emptyset$.) Therefore, $b\in B(b_i,\lambda/2)$ for some $i\in \{1,\cdots,N\}$ where $b_i\in A\cap B$.  Now, $|v-b_i|\le |x| +\lambda/2 <\lambda$, and so $v\in \{A\cap B\}^{\lambda}$.


\vskip.125in 

\newpage

\begin{thebibliography}{8}








\bibitem{EIT11} S. Eswarathasan, A. Iosevich and K. Taylor, {\it Fourier integral operators, fractal sets and the regular value theorem}, Adv. Math. \textbf{228} (2011), 2385-2402. 


\bibitem{Fal86} K. J. Falconer, {\it On the Hausdorff dimensions of distance sets} Mathematika \textbf{32} (1986), 206-212.

\bibitem{Falc86} K. J. Falconer, {\it The geometry of fractal sets}, Cambridge Tracts in Mathematics, \textbf{85} Cambridge University Press, Cambridge, (1986). 

\bibitem{Falc94} K. J. Falconer, {\it Sets with large intersections}, J. London Math. Soc. \textbf{49} (1994), 267-280.

\bibitem{Fe69} H. Federer, {\it Geometric Measure Theory}, 2.10.29, (1969).










\bibitem{IJL10} A. Iosevich, H. Jorati and I. Laba, {\it Geometric incidence theorems via Fourier analysis}, Trans. Amer. Math. Soc. \textbf{361} (2009), 6595-6611.








\bibitem{Mat84} P. Mattila, {\it Hausdorff dimension and capacities of intersections of sets in $n$-space}, Acta Math. 
\textbf{152} (1984), 77-105.

\bibitem{Mat85} P. Mattila, {\it On the Hausdorff dimension and capacities of intersections}, Mathematika \textbf{32} (1985), 213-217.

\bibitem{Mat87} P. Mattila {\it Spherical averages of Fourier transforms of measures with finite energy: dimensions of intersections and distance sets}, Mathematika \textbf{34} (1987), 207-228.

\bibitem{M95} P. Mattila, {\it Geometry of sets and measures in Euclidean spaces}, Cambridge University Press, \text{volume} 44, (1995). 












\end{thebibliography}
\end{document}